\input amstex
\documentstyle{amsppt}
\voffset-,6in
\magnification=\magstep1
\NoBlackBoxes
\topmatter
\title On a Certain Class of $g$-Functions for Subshifts
\endtitle
\author  Wolfgang Krieger
\endauthor
\affil 
Institute for Applied Mathematics\\
University of Heidelberg\\
Im Neuenheimer Feld 294
69120 Heidelberg,
Germany\\
\endaffil
\abstract
A property $(D)$ of subshifts was defined in: Wolfgang Krieger, On
$g$-functions for subshifts, IMS Lecture Notes- Monograph Series, Vol. 48,
Dynamics \& Stochastics (2006) 306 - 316, arXiv:math.DS/0608259. With a view
towards
a theory of $g$-fnctions beyond the case of finte type subshifts partially
defined continuous $g$-functions of property $(D)$ subshifts are studied.
\endabstract
\endtopmatter





\bigskip
\heading 1.Introduction
\endheading

Let $\Sigma $ be a finite alphabet, and let $S$ denote the shift on
$\Sigma
^{\Bbb Z}$,
$$
S((x_{i}) _{i \in \Bbb Z}) = ((x_{i+1}) _{i \in \Bbb Z}),\qquad (x_{i}) _{i
\in \Bbb
Z} \in \Bbb Z.
$$
A closed $S$-invariant set $X \subset  \Sigma ^{\Bbb Z} $ with the restriction
of $S$ acting on it, is called a subshift. A finite word is said to
be admissible for a  subshift if it appears in a point of the subshift. A
subshift is uniquely determined by its set of admissible words. A subshift is
said to be of finite type if its admissible words are defined by
excluding finitely many words from appearing as subwords in them. Subshifts are
studied in symbolic dynamics. For an
introduction to symbolic dynamics see \cite {Ki} and \cite {LM}.

$g$-functions and $g$-measures were introduced (under this name) by Mike Keane
\cite {Ke1}, \cite {Ke2}. For other contexts compare e.g. \cite {DF}\cite
{H}\cite {K}\cite {OM}. A
substantial theory of
$g$-functions and
$g$-measures of subshifts of finite type exists.
$g$-functions for more general subshifts were considered in \cite {Kr2} (see
also \cite {Mai}, p. 24). Here we continue this line of investigation by
introducing and examining a class of partially defined continuous
$g$-functions for property $(D)$ subshifts.

Before formulating the relevant notions we recall the notation that we use for
subshifts. Given a subshift $X \subset \Sigma^{\Bbb Z}$
we set
$$
	x_{[i,k]} = (x_{j})_{i \leq j \leq k}, \quad  		x \in X,  
			i,k\in \Bbb Z ,  i \leq k,
$$
and
$$
X_{[i,k]} = \{ x_{[i,k]} : x\in X \} . 
$$
and use similar notation also if indices range in semi-infinite intervals.
Blocks can also 
stand for the words
they carry.
We denote
$$
\Gamma ^{+}_{n} (x^{-} ) = \{ b \in X_{[1, n ]} : (  x^{-} , b) \in  X_{(-
\infty, n ]}  \}, \qquad n \in \Bbb N,
$$
$$
 \Gamma ^{+} (x^{-} ) =  \bigcup_{  n \in \Bbb N  } \Gamma ^{+}_{n} (x^{-} ),
$$
$$
\Gamma ^{+}_{\infty} (x^{-} ) = \{ x^{+} \in X_{[1, \infty )} : (  x^{-} ,
x^{+}) \in  X  \}   , \qquad  
x^{-} \in  X_{(- \infty, 0 ]}
$$
and
$$
\Gamma ^{+}_{n} (a ) = \{ b \in X_{[1, n ]} : (  a , b) \in  X_{(- k, n ]} 
\}, \qquad n \in \Bbb N,
$$
$$
 \Gamma ^{+} (a ) =  \bigcup_{  n \in \Bbb N  } \Gamma ^{+}_{n} (a),
$$ 
$$
\Gamma ^{+}_{\infty} (a ) = \{ x^{+} \in X_{[1, \infty )} : (  a , x^{+}) \in
 X  \}   , \qquad  
a \in  X_{(- k, 0 ]}, k \in \Bbb Z_{+}.
$$ 
$\Gamma^{-}$ has the time symmetric meaning. We denote
$$
\omega^{+}_{n}(a) = \bigcap_{x^{-} \in \Gamma^{-}_{\infty}(a)}\{ b \in x_{[1,
n]}: ( x^{-},  a,  b) \in X_{( - \infty, n]}\},
$$
$$
\omega_{n} ^{+} (x^{-} ) =  \bigcup_{  k \in \Bbb N  } \omega^{+}_{n}
(x^{-}_{(-k. 0]}),\qquad  
x^{-}  \in  X_{(- \infty, 0 ]}, n \in \Bbb N.
$$

In \cite{Kr2} a
property $(D)$
of subshifts was introduced. 
A subshift $X \subset \Sigma ^{ \Bbb Z}$ has property $(D)$ if for
all
admissible words $b\sigma$ of $X$ there exists a word $a \in 
\Gamma^{-}(b)$  
such that $\sigma \in \omega ^{+}_{1} (ab)$. Property $(D)$ implies that for for
all
admissible words $b$ and all words $c\in \Gamma^{+}_{n}(b), n \in \Bbb
N$, there exists a word $a \in 
\Gamma^{-}(b)$ 
such that $c \in \omega ^{+}_{n} (ab)$ (\cite {Kr2}, Lemma 4.3.). Property
$(D)$ is an invariant of
topological conjugacy (\cite {Kr2}, Proposition 4.3). There exist coded systems
\cite {BH} without property $(D)$.

Let 
$X \subset \Sigma^{\Bbb Z}$
be a subshift. With residual subsets $D^{-}$ of 
$X_{(-\infty, 0]}$
we consider continuous functions
$$
g : D^{-}   \times \Sigma \to [0,1], 
$$
such that
$$
 \sum_{\sigma \in \Gamma^ {+}_{1}(x^{-})}   g(x^{-}   , \sigma) = 1,\qquad
x^{-} \in  D^{-}  . 
$$
Here we always assume that there is no point in
$X_{(-\infty, 0]}-D^{-}$
to which the mapping
$x^{-} \to (g(x^{-}   , \sigma))_{\sigma \in \Sigma} \ (x^{-} \in  D^{-})$
can be extended by continuity. $D^{-}$ is then a $G_{\delta}$ that we denote by
$ D^{-}(g).$
We say that $g$ is a residually defined $g$-function of $X$ if 
$$
X_{(- \infty, 1]} = \overline { \{(x^{-}, \alpha): x^{-}\in  D^{-}(g), \alpha
\in \Gamma^{+}_{1}(x^{-}), g(x^{-} , \alpha) > 0 \}}.
$$
Note that a 
residually defined $g$-function of $X$ is uniquely determined by its
restriction to any dense subset of $D^{-}(g) \times \Sigma$. In section 2
we show that a subshift has a residually defined $g$-function if and only if
it has property $(D)$.

Given a residually defined $g$-function of a property $(D)$ subshift
$X\subset \Sigma^{\Bbb Z}$,
we say that an invariant probability measure  $\mu$ $X$ is a $g$-measure 
for $g$, if
$$
 \mu (\{ x \in X: x_{( \infty, 0]} \in D^{-}(g) \})= 1, \tag 1
$$
 and if
$$
\align
\mu(\{ &x \in X: x_{( \infty, 0]} \in D^{-}(g), x_{[-k,1]}= (a, \alpha)\}) =
\\ 
&\int_{\{x\in X: x_{( \infty, 0]} \in D^{-}(g), x_{[-k,0]}= a\}}g(x_{(-\infty,
0]}, \alpha) d\mu, 
  a \in X _{[-k,0]} , \alpha \in \Gamma^{+}_{1}(a), k \in \Bbb N. \tag 2
\endalign
$$ 
For a residually defined $g$-function of a property $(D)$ subshift $X \subset
\Sigma^{\Bbb Z}$ we set
$$\align
 E(g) = \bigcap_{i \in \Bbb Z} (\{x\in X&:x_{(-\infty, i]} \in
S^{-i}D^{-}(g),g(x_{(-\infty, i]}, x_{i+1}) > 0 \}\cap 
\\
& \bigcap_{ k \in \Bbb N} \{x \in
X:
\{(x_{(-\infty, i]},
a): a \in \Gamma^{+}_{k}(x_{(-\infty, i]})\} \subset S^{-i - k}D^{-}(g)\}).
\endalign$$
We will see in Section 3 that $(g)$ is a dense $G_{\delta}$. Note that (1) and
(2) imply that for a
$g$-mesaure $\mu$ of a residually defined $g$-function $g$, $ \mu(E(g)) = 1$.

We recall some terminology for labeled directed graphs. A directed graph
with vertex set $\Cal V$ and edges carrying labels taken
from a finite alphabet $\Sigma$ is called a Shannon graph if the labeling is
1-right resolving in the sense that for all $ V \in  \Cal V$ and $ \sigma
\in \Sigma$ there is at most one edge leaving $ V $ that carries the label
$\sigma$.
Denote the set of vertices that have an outgoing edge that carries the label
$\sigma$ by $\Cal V(\sigma)$, and for $V \in \Cal V(\sigma) $ denote by
$\tau(\sigma)V$
the final vertex of the edge that leaves $V$ and carries the label $\sigma$.
We call the mappings $V \to \tau(\sigma)V \ (v \in  \Cal V (\sigma)), \sigma \in
\Sigma$ the transition rules of the Shannon graph.
For a word $a \in \Sigma^{n}, n \in \Bbb N$, and for $V \in  \Cal V, 
\tau(a)V$
will denote the end vertex of the path that leaves $V$ and carries the label
sequence $a$ (if such a path exists).
Call a set $\Cal W \subset \Cal
V $ transition complete if for $\sigma \in \Sigma, V \in \Cal 
W 
\cap \Cal V (\sigma)$ implies that also $\tau(\sigma)W \in \Cal W$. Every
transition complete set $\Cal W \subset \Cal V $
determines
a Shannon
graph with transition rules that are inherited from the Shannon graph 
$\Cal V
$. 
For a vertex $V$ of a Shannon graph its forward context is defined as the set
of label sequences of the paths that leave $V$. A Shannon graph is called
forward separated if distinct vertices have distinct forward contexts.  A
Shannon graph is
said to present a subshift $ X \subset  \Sigma ^{\Bbb
Z} $ if the
set of admissible words of the subshift coincides with the set of label
sequences of finite paths in the graph.
A compact Shannon graph is a Shannon graph whose vertex set carries a compact
topology such that the sets $\Cal V(\sigma)$ are open, and such that the
transition rules are
continuous. 
To a compact Shannon graph $\Cal G$ such that every vertex has an incoming
an an outgoing edge and such that the sets $\Cal V (\sigma ), \sigma \in  
\Sigma $ are compact,
there is associated a topological Markov chain
$$
M(\Cal G) = \bigcap_{i\in \Bbb Z} \{ (x_{i}, V_{i})_{i \in \Bbb Z} \in (\Sigma
\times \Cal V )^{\Bbb Z}: V_{i+1} = \tau(x_{i+1})V_{i}  \}.
$$
This topological Markov chain $ M(\Cal G) $ projects onto the subshift that
contains the label sequences of the two-sided infinite paths on $\Cal G$. 

For a finite alphabet $\Sigma$ denote by $\Cal V(\Sigma)$ the set of closed
subsets of $\Sigma ^{\Bbb N}$ with its compact Hausdorff subset topology.
$\Cal V (\Sigma)$ is the vertex set of a compact Shannon graph $\Cal G
(\Sigma)$. The set $\
\Cal V (\sigma)$ in $\Cal V(\Sigma)$ is the set of $V\in\Cal V (\Sigma)$ that
contains a sequence that
starts with $ \sigma \in \Sigma$ and one has the transition rules
$$
\tau(\sigma)V = \{ v_{[1, \infty)}: v \in V, v_{1} = \sigma \}, \qquad V \in 
\Cal V(\sigma), \sigma \in \Sigma.
$$
A transition complete  sub-Shannon graph of $\Cal G
(\Sigma)$ is forward separated and for every forward separated  Shannon graph
with vertex set $\Cal V$ one has the isomorphic transition complete
sub-Shannon graph of $\Cal G
(\Sigma)$ with vertex set $\{\Gamma^{+}(V): V \in \Cal V\}$.

For a finite alphabet $\Sigma$ denote by $\Cal M (\Sigma )$ the set of
probability measures on $\Sigma^{ \Bbb N }$ 
 with its weak *-topology. With the notation
$$
C(a) = \bigcap_{ 1 \leq i \leq n } \{ (x_{i}) _{i \in \Bbb N} \in \Sigma^{
\Bbb N }: a_{i}= x_{i}\}, \qquad (a_{i})_{1 \leq i \leq n} \in \Sigma^{N}, n
\in \ \Bbb N, 
$$
$$
\Cal M (\Sigma )(\sigma) = \{ \mu \in \Cal M (\Sigma) : \mu (C(\sigma)) > 0 \},
\qquad \sigma \in \Sigma,
$$
let for  $ \mu \in \Cal M (\Sigma )(\sigma)$, $\tau(\sigma)\mu$ be
equal to the conditional measure of $\mu$ given $C(\sigma)$,
$$
\tau(\sigma)\mu (C(b)) = \frac{\mu(C(\sigma, b))}{ \mu(C(\sigma))} , \qquad b
\in \Sigma^{N}, N \in \Bbb N.
$$
In this way $\Cal M(\Sigma )$ has been turned into a compact Shannon graph
with the
transition rules $(\tau(\sigma))_{\sigma \in \Sigma}$.


In Section 2 we construct for a property $(D)$ subshift $X$ an invariantly
associated presenting compact Shannon graph $\Cal G_{D}$ whose vertex set is a
transition  complete subset of $\Cal V (\Sigma)$.

In Section 3 we introduce a notion of
topological conjugacy for residually
defined $g$-functions and we list some invariants of this conjugacy. We
show that a topological conjugacy of $g$-functions carries $g$-measures into
$g$-measures.

In Section 4 we describe a one-to-one correspondence between residually
defined $g$-functions of property $(D)$ subshifts and certain compact
Shannon graphs whose vertex sets are transition complete subsets of $\Cal M
(\Sigma)$ and that we call
residually contractive Shannon graphs. 

\heading 2. A Shannon graph for property $(D)$ subshifts
\endheading

Set for a subshift
$X\subset \Sigma^{\Bbb Z}$,
$$
E(X) = \bigcap_{i\in \Bbb Z} \{ x \in X: x_{i} \in  \omega_{1}^{+}(x)_{(-\infty,
i)})\}.
$$
\proclaim{Proposition 2.1}
Let 
$X\subset \Sigma^{\Bbb Z}$
be a subshift with property $(D)$. Then E(X) is residual in X.
\endproclaim
\demo{Proof}
Setting
$$
E_{i,k}(X) =
\bigcup_{\{a \in X_{[i-k.i]}: a_{i} \in \omega_{1}^{+}(a_{[i-k.i)})\}} \{x \in
X :x_{[i-k.i]}   = a\},\qquad  i \in \Bbb Z, k \in \Bbb N.
$$
one can write
$$
E(X) = \bigcap_{i \in \Bbb Z} \bigcup_{k \in \Bbb N}  E_{i,k}(X) .
$$
This shows that $E(X)$ is a
$G_{\delta}$.
We prove that the set
$\bigcup_{k \in \Bbb N}  E_{0,k}(X)$ 
is dense in $X$. For this, let
$n \in \Bbb N, c \in X_{[-n, n]}$.
By property $(D)$ of $X$ there is an $ m \in \Bbb N$ and a 
$b \in \Gamma^{-}_{m}(c)$
such that 
$c \in \omega_{n}^{+}(b).$
Then
$$
\{ x \in X: x_{[-n-m,n]} = bc\} \subset   E_{0,n+m}(X)  .\qed
$$
\enddemo
We set
$$
\Omega^{+}_{\infty}(x^{-})= \bigcap_{n \in \Bbb N} \{x^{+} \in X_{[1,
\infty)}: x^{+}_{[1,n]} \in
\omega^{+}_{n}(x^{-})\},\qquad x^{-} \in X_{(-\infty, 0]}.
$$
\proclaim{Lemma 2.2} 
$\Omega^{+}_{\infty}(x^{-})_{[1,n]} =\bigcap_{m \geq n}
\omega^{+}_{m}(x^{-})_{[1,n]}, \qquad n \in \Bbb N, x^{-} \in X_{(-\infty,
0]}.   $
\endproclaim
\demo{Proof}
For the proof, we let 
 $n \in \Bbb N, x^{-} \in X_{(-\infty,
0]}, $ and
$$
a \in \bigcap_{m \geq n} \omega^{+}_{m}(x^{-})_{[1, n]},
$$
and we produce an 
$x^{+} \in \Omega^{+}_{\infty}(x^{-})$
such that
$x^{+}_{[1,n]} = a$
in the following way: It is
$$
\bigcap_{m \geq k} \omega^{+}_{m}(x^{-})_{[1,k]} = ( \bigcap_{m > k}
\omega^{+}_{m}(x^{-})_{[1,k+1]})_{[1,k]},\qquad  k \in \Bbb N .
$$
We can therefore choose inductively a sequence
$a^{(k)}, k \geq n,a^{(n)} = a,$
$$
 a^{(k)}  \in \bigcap_{m \geq k}
\omega^{+}_{m}(x^{-})_{[1,k]},\qquad k \geq n,
$$
such that
$$
a^{(k+1)}_{[1, k]} = a^{(k)},\qquad k \geq n.
$$
Let
$x^{+} \in \Omega^{+}_{\infty} (x^{-})$
be given by
$ x^{+}_{[1,k]}= a^{(k)}, k \geq n.$
\qed
\enddemo
\proclaim{Lemma 2.3}
The sets
$\Omega^{+}_{\infty} (x^{-}), x^{-} \in X_{(- \infty, 0]},$
are closed.
 \endproclaim
\demo{Proof}
Apply Lemma 2.2.
\qed
\enddemo
\proclaim{Lemma 2.4}
 $$
\{x^{-} \in X_{(-\infty,0]}: \Omega^{+}_{\infty} (x^{-} )\neq \emptyset\} =
\bigcap_{n \in \Bbb N} \{x^{-} \in X_{(-\infty,0]} : \omega^{+}_{n} (x^{-} )
\neq \emptyset  \}
$$
\endproclaim
\demo{Proof}
 Apply Lemma 2.2.
\qed
\enddemo

We set
$$
E^{-}(X) =
\{x^{-} \in X_{(-\infty, 0]}: \Omega^{+}_{\infty}(x^{-}) \neq \emptyset\} \cap
\bigcap_{i \in \Bbb Z_{-}} \{  x^{-} \in X_{(-\infty, 0]}: x_{i} \in
\omega^{+}_{1}(x_{(-\infty, i)})\}.
$$

\proclaim{Proposition 2.5}
Let 
$X \subset \Sigma^{\Bbb Z}$
be a subshift with property $(D)$. Then $E^{-}(X)$ is residual in $X_{(-\infty,
0]}.$
 \endproclaim
\demo{Proof}
One has
$$
\{x^{-} \in X_{(-\infty, 0]}: \omega^{+}_{n}(x^{-}) \neq \emptyset\} =
\bigcup_{k \in \Bbb N} \bigcup_{a\in X_{[1,n]} } \{x^{-} \in  X_{(- \infty, 0]}:
a \in
\omega^{+}_{n}(x^{-}_{-k, 0]})\}, n \in \Bbb N,
$$
and it is seen from Lemma 2.4 that the set 
$\{x^{-} \in X_{(-\infty, 0]}: \Omega^{+}_{\infty}(x^{-}) \neq \emptyset\}$
is a $G_{\delta}$.

To see that $E^{-}(X)$ is dense in 
$X_{(-\infty, 0]},$ 
observe that $E^{-}(X) $ is the image of $E(X)$ under the projection  of 
$X$ onto $X_{(-\infty, 0]}$ : If $x \in E(X)$ then 
$x_{[1, \infty)} \in  \Omega^{+}_{\infty}(x^{-})$, and if 
$x^{-} \in E^{-}(X)$
and
$x^{+} \in  \Omega^{+}_{\infty}(x^{-})$
then
$( x^{-} , x^{+} ) \in E(X)$.
Apply Proposition 2.1.
\qed
\enddemo

\proclaim{Corollary 2.6}
 A subshift
$X \subset \Sigma^{\Bbb Z}$
has property $D$ if and only if $E(X)$ is dense in $X$.
\endproclaim
\demo{Proof}
Assume that $E(X)$ is dense in $X$ and let $a  \in X_{[-k,0]}, k \in \Bbb 
N.$
Let $x\in E(X), a = x_{[-k,0]}$. Then there is an $l \in \Bbb N$ such that
$x_{0} \in \omega_{1}^{+}(x_{[-k-l,0)}    ),$ and it is seen that $X$ has
property
$(D)$. 
\qed
\enddemo
We note that similarly a subshift
$X \subset \Sigma^{\Bbb Z}$ has property $(D)$ if and only if
$E^{-}(X)$
is dense in $X_{(-\infty, 0]}$.

One obtains for a property $(D)$ subshift $X \subset \Sigma^{\Bbb Z}$ a
presenting compact
Shannon graph $\Cal G_{D}(X)$ whose vertex set is  given the closure of the set 
$\{\Omega^{+}_{\infty}(x^{-}): x^{-} \in E^{-} (X) \}$. By means of Lemma 2.1
of \cite {Kr1} it can be proved that the compact Shannon graph $\Cal G_{D}(X)$
is invariantly associated to the subshift $X$. By this is meant that a
topological conjugacy of a property $(D)$ subshift $X \subset \Sigma^{\Bbb Z}$
onto a subshift $\widetilde{X}\subset \widetilde{\Sigma}^{\Bbb Z}$ induces a
topological
conjugacy of the topological Markov chain $M(\Cal G_{D}(X))$ onto the
topological Markov chain $M(\Cal G_{D}(\widetilde{X}))$. In \cite {Mat}
Matsumoto has introduced the notion of a $\lambda$-graph system and the notion
of the strong equivalence of $\lambda$-graph systems. The compact
Shannon graphs 
$\Cal G_{D}(X)$, that are associated to property $(D)$ subshifts $X$ can be
replaced by equivalent structures that are forward separated Shannon
$\lambda$-graph system \cite{KM}. In this way topologically conjugate property
$(D)$
subshifts give rise to strong shift equivalent $\lambda$-graph
systems. 

\heading 3. Residually defined $g$-functions
\endheading

\proclaim{Proposition 3.1}
A subshift admits a residually defined g-function if and only if it has
property (D).
\endproclaim
\demo{Proof}
Let $g$ be a residually defined $g$-function of the subshift $X\subset
\Sigma^{\Bbb Z}.$ Also let  $(a, \alpha) \in X_{(-k,1]}, k \in \Bbb N$. There
exists
a point $x\in E(g) \cap C(a, \alpha)$. $g(x_{(-\infty, 0]}, \alpha)
> 0$ implies by Lemma 2.2 of \cite {Kr2} that $\alpha \in
\omega^{+}(x_{(-\infty, 0]}$
and 
therefore there is a $b \in \Gamma^{-}(a)$ such that $\alpha \in
\omega_{1}^{+}(ba)$,
which means that $X$ has property $(D)$.

The converse follows from Lemma 2.3 of \cite {Kr2}. 
\qed
\enddemo
\proclaim {Lemma 3.2}
Let $ X \subset \Sigma^{\Bbb Z}$ be a subshift, let $k \in \Bbb N    $ , and
let $A^{-}$ be a residual subset of $X_{(-\infty, k]}$. Then the set
$$
\{ x^{-} \in X_{(-\infty, 0]}: \{(  x^{-} , a): a \in \Gamma^{+}_{k}( x^{-}
)\} \subset S^{-k} A^{-} \}
$$
is residual in $ X_{(-\infty, 0]}  $.
\endproclaim
\demo{Proof}
For 
$a \in X_{[1,k]}$ the set 
$  S^{-k} A^{-}\cap \{ (  x^{-} , a):  x^{-}\in \Gamma^{-}_{\infty}(a)\}$ is
residual
in
$\{ (  x^{-} , a): x^{-}  \in \Gamma^{-}_{\infty}(a)\} $, and this implies
that the set $ \{ x^{-}  \in \Gamma^{-}_{\infty}(a): (  x^{-} , a) \in 
S^{-k}A^{-}\}
$
is
residual in $\Gamma^{-}_{\infty}(a)$. Consequently, the set
$$
\bigcup_{a \in   X_{[1,k]}  }\{ x^{-}  \in \Gamma^{-}_{\infty}(a):   (  x^{-}
, a)   \notin  S^{-k}A^{-} \}
$$
is a countable union of nowhere dense closed subsets of  $X_{(-\infty, 0]} $,
and the lemma follows.
\qed
\enddemo

Note that, as a consequence of Lemma 3.2, $E(X)$ is a dense $ G_{\delta}$.
Note also that by Lemma
2.1 of \cite {Kr2} $ E(g) \subset E(X)$. We
set
$$
F_{i} (x) = \{ y \in X: y_{(-\infty, i]} = x_{(-\infty, i]} \}, \qquad x \in
E(g),i \in \Bbb Z,
$$
and 
$$
F_{i}(x,a) = \{y \in F_{i}(x): y_{(i, i+n]}= a\},
\qquad
a \in \Gamma^{+}
_{n}(x_{(-\infty, i]}), n \in \Bbb N, x \in E(g).
$$
We define a probability measure
$\mu^{(g)}_{i}(x_{(-\infty, i]}) )$ on $  F_{i} (x)  $ by setting
$$\align
\mu^{(g)}_{i}(x_{(-\infty, i]})(F_{i}(x,a))    & =
\prod_{1 \leq k \leq N} g(( x_{(-\infty, i]} , a_{[i+1, i+k)}), a_{i+k}),\\
a &\in \Gamma^{+}_{N}(x_{(-\infty, i]}), N \in \Bbb N, x \in E(g).
\endalign
$$
By (\cite {Kr2} Lemma 2.1),
$$
\mu^{(g)}_{i}(x_{(-\infty, i]})
)(\Omega^{+}_{(\infty}(\mu^{(\varphi)}_{i}(x_{(-\infty, i]}) ))) = 1, \qquad
x\in E(g).
$$

A measure theoretic structure is now apparent that we describe. Let $(F, \Cal
F)$ be a measurable space, and let $F_{i} \in  \Cal F, i \in \Bbb Z$, be such
that 
$$
 F_{i} \supset   F_{i+1    },\qquad i \in \Bbb Z,
$$
and such that
$$
F = \bigcup_{k \in  \Bbb N}F_{-k}.
$$
Denote by $\Cal H ((F_{i}  )_{ i \in \Bbb Z  })$ the set of $\sigma$-finite
measures on $(F, \Cal F)$ such that
$$
0 < \mu(F_{i} ) < \infty,\qquad i \in \Bbb Z.
$$
Note that
$$
\Cal H (F_{i}  )_{ i \in \Bbb Z  })  = \Cal H ((F_{jL}  )_{ j \in \Bbb Z  }),
\qquad  L \in \Bbb N,
$$
and that
one has for 
$\widetilde{F}_{i} \in \Cal F,   i \in \Bbb Z $,
 such that
$$  F_{i+1} \supset    \widetilde{F}_{i}       \supset F_{i}, \qquad    i \in
\Bbb Z ,
$$
that
$$
\Cal H ((F_{i}  )_{ i \in \Bbb Z  }) = \Cal H ((\widetilde{F}_{i} )_{ i \in
\Bbb Z  }) .
$$
Call $ \mu, \mu^{\prime} \in \Cal H ((F_{i}  )_{ i \in \Bbb Z  })$ scale
equivalent, $ \mu \sim_{sc}\mu^{\prime}$, if for some $ \alpha > 0, 
\mu^{\prime}   = \alpha  \mu .$ Denote by $ H ((F_{i}  )_{ i \in \Bbb Z  })  $
the
set of sequences 
$(\mu_{i})_{i \in \Bbb Z },$ where $  \mu_{i}  $ is a probability measure on $
F_{i} $, such that
$$
   \mu_{i} ( F_{i+1}  )   > 0 ,\qquad i \in \Bbb Z ,
$$
and such that
$$
   \mu_{i+1} = \frac {  \mu_{i}  \restriction  F_{i+1}}{ \mu_{i} ( F_{i+1})}  ,
   \qquad i \in \Bbb Z .
$$
There exists a bijection $\eta$ of $[  \Cal H ((F_{i}  )_{ i \in \Bbb Z  }) 
]_{\sim_{sc}}$ onto $  H ((F_{i}  )_{ i \in \Bbb Z  })  $ that is given by
$$
\eta([ \mu  ]_{\sim_{sc}}) = ( \frac{ \mu  \restriction  F_{i}  }{  \mu  (
F_{i}) })_{ i \in \Bbb Z  },\qquad \mu
\in   \Cal H ((F_{i}  )_{ i \in \Bbb Z  }).
$$
For an element $(\mu_{i})_{i \in \Bbb Z }$ of $  H ((F_{i}  )_{ i \in \Bbb Z
 })  $ an element of the scale equivalence class that is its inverse under the
mapping $\eta$ is, for instance, the measure $ \mu
\in   \Cal H ((F_{i}  )_{ i \in \Bbb Z  }$ that is given by setting
$$
 \mu  \restriction  F_{0} = \mu_{0},
$$
and
$$
\mu  \restriction (F_{-k}-F_{-k+1}) = \frac { \mu_{-k}  \restriction
(F_{-k}-F_{-k+1})  }{ \mu_{-k} (F_{0} )  }, \qquad k \in \Bbb N.
$$

For a residually defined $g$-function $g$
of the property $(D)$ subshift $X \subset \Sigma^{\Bbb Z}$ one has
$$
( \mu^{(g)}_{i}(x_{(-\infty,i]}))_{i\in \Bbb Z} \in H(( F_{i} (x))_{i\in \Bbb Z}
), \qquad x \in E(x).
$$
We say that a residually defined $g$-function $\widetilde{g}$ of a property
$(D)$ subshift $\widetilde{g}$
is is the image of a residually defined 
$g$-function $g$ of a property $(D)$ subshift 
$X \subset \Sigma^{\Bbb Z}$ under a 
topological conjugacy $\varphi: X  \to \widetilde{X }  $ if
$$
  \varphi  (E(g)) = E(\widetilde{g}),
$$
and if
$$
\eta^{-1}((( \mu^{(g)}_{i}(x_{(-\infty,i]}))_{i\in \Bbb Z} ) =
\eta^{-1}(((  \varphi^{-1}\mu^{(\widetilde{g})}_{i}(
\varphi(x)_{(-\infty,i]}))_{i\in \Bbb Z} ), \qquad  x\in E(g). 
$$
If for some $L \in \Bbb N$ the topological conjugacy $\varphi^{-1} $ has the
coding window $[-L, L]$
then (1) holds if and only if
$$
\varphi^{-1}
 ( \mu^{(\widetilde{g})}_{0}(\varphi(x)_{(-\infty,0]}))=
 \frac 
{\mu^{(g)}_{-L}(x_{(-\infty,-L]})\restriction\varphi^{-1}(F_{0}(\varphi(x)))}
{ \mu^{(g)}_{-L}(x_{(-\infty,-L]})(\varphi^{-1}(F_{0}(\varphi(x))))}, \qquad
x \in E(g). \tag 2
$$
For the $g$-function $\widetilde{g}$ itself it follows that
$$
\widetilde{g}(  \widetilde{x}_{(-\infty,0]}, \alpha) = \frac 
{\mu^{(g)}_{-L}(x_{(-\infty,-L]})(\varphi^{-1}(F_{0}(\varphi(x), \alpha)))}
{ \sum_{\beta \in
\Gamma^{+}_{1}(x_{(-\infty,0]} )  }
\mu^{(g)}_{-L}(x_{(-\infty,-L]})(\varphi^{-1}(F_{0}(\varphi(x), 
\beta)))}, \tag 3
$$
$$
 \qquad \alpha \in \Gamma^{+}_{1}(x_{(-\infty,0]}, x \in E(g). 
$$
\proclaim{Theorem 3.3}
The finiteness of the range of a residually defined $g$-function $g$ of a
property $(D)$ subshift is an invariant of topological conjugacy.
\endproclaim
\demo{Proof}
Inspect formula (3)
\qed
\enddemo

\proclaim {Theorem 3.4}
For residually defined $g$-functions $g$ of  property $(D)$ subshifts the
finiteness and countable infinity of the range of the mapping
$$
x_{(-\infty,0]} \to \mu^{(g)}_{0} (x_{(-\infty,0]})  ( x \in E(g))
$$
are invariants of topological conjugacy.
\endproclaim
\demo{Proof}
Inspect formula (2).
\qed
\enddemo

We recall the notion of a bipartite subshift. Let $\Delta $ and $\widetilde
{\Delta}$ be finite disjoint alphabets, and let $Y \subset ( \Delta   \cup
\widetilde {\Delta} )^{\Bbb Z}$ be a subshift. $Y$ is called bipartite if the
admissible words of length two of $Y$ are contained in $\Delta  \widetilde
{\Delta} \cup \widetilde {\Delta} \Delta $. If $Y$ is bipartite then 
$S_{Y}^{2}$
leaves the sets 
$$
X = \{ (y_{i})_{i\in \Bbb Z} \in Y: y_{0} \in  \Delta     \},
$$
and
$$
\widetilde{X }= \{ (y_{i})_{i\in \Bbb Z} \in Y: y_{0} \in  \widetilde {\Delta}
   \},
$$
invariant. Let $S$ resp $\widetilde{S}$ denote the restriction of  
$S_{Y}^{2}$  to $X$
resp. to $ \widetilde {X}$. $(X,S)$ and  $(\widetilde {X},\widetilde {S})$ are
topologically conjugate: a
topological conjugacy of $X$ onto $\widetilde{X}$ is given by the restriction of
$S_{Y}$ to
$X$. Denote the set of words in $  \Delta  \widetilde {\Delta}   $  resp. in $
 \widetilde {\Delta} \Delta $ that are admissible for $Y$ by $\Sigma$ resp. by $
 \widetilde{ \Sigma}     $.
One has $X \subset \Sigma ^{\Bbb Z},    
\widetilde{X} \subset \widetilde { \Sigma} ^{\Bbb Z}$, and one has the
injections
$$
\psi : \Sigma \hookrightarrow \Delta \widetilde{\Delta}, \quad \widetilde
\psi: 
\widetilde \Sigma \hookrightarrow \widetilde\Delta \Delta.
$$
By applying $\psi$ and $\tilde \psi $ symbol by symbol one extends their
domain of 
definition to finite words and right-infinite sequences.
$\psi$ and $\widetilde \psi $ satisfy the relation
$$
\iota^{-}( \iota^{+}    ( \psi( X_{[1,2]}))) = \widetilde{\psi}
(\widetilde 
\Sigma), \tag 2
$$
and are called specifications. 
Conversely, let $X \subset \Sigma^{\Bbb Z}$ and $ \widetilde{ X} \subset
\widetilde{\Sigma}^{\Bbb Z}   $ be subshifts, and let $\Delta$ and $\widetilde{
\Delta}     $ be disjoint finite alphabets and let there be given
injections
$$
\psi : \Sigma \hookrightarrow \Delta \widetilde{\Delta}, \quad 
\widetilde \psi :\widetilde \Sigma \hookrightarrow
\widetilde\Delta\Delta , 
$$
that are specifications, that is, they satisfy (2). Then
$$
\psi (X) \cup   \widetilde{\psi}(\widetilde{X})    \subset \{ \Delta 
\widetilde{\Delta }   \cup 
\widetilde{\Delta } \Delta  \}^{\Bbb Z}
$$
is a bipartite subshift and $X$ and   $ \widetilde{ X} $ are topologically
conjugate. (1) implies that also 
$$
\iota^{-}(\iota^{+}(\psi (\widetilde X_{[1,2]}))) = \psi (\Sigma)
$$
and one has a 2-block code given by
$$
a \to \widetilde{\psi}^{-1}(\iota^{-}(\iota^{+}(\psi (a)))), \qquad a
\in X_{[1,2]},
$$
that implements a topological conjugacy of $X$ onto $ \tilde{X}$ with the
inverse
given by the 2-block map
$$  
\widetilde{a }\to \psi^{-1}(\iota^{-}(\iota^{+}(\widetilde{\psi} 
(\widetilde{a })))), 
\qquad\widetilde{a} \in\widetilde{
X}_{[1,2]}.
$$
Topological conjugacies that arise in this way are called bipartite codings. If
a bipartite coding exists between two subshift presentations then these are
called bipartitely related. 
According to a theorem of Nasu \cite{N} subshifts $X \subset \Sigma 
^{\Bbb Z}     $ and $\widetilde{X} \subset \widetilde{\Sigma} 
^{\Bbb Z}     $ are topologically conjugate 
if and only if there is a chain $ X(k) \subset 
\Sigma(k) 
^{\Bbb Z}     , 0 \leq k \leq 
K, K \in \Bbb N$, of subshifts,    $ X(0) = X,  X(K) = \widetilde{X}  $ such
that 
 $ X(k) \subset 
\Sigma(k) 
^{\Bbb Z}     $ and $  X(k+1) \subset 
\Sigma(k+1) 
^{\Bbb Z}     $ are bipartitely related, $ 0 \leq k <
K $.

\proclaim{Lemma 3.5}
Let $g$ be a residually defined $g$-function of the property (D) subshift $X
\subset \Sigma^{\Bbb Z}$, and let $\widetilde{g}$ be a residually defined
$g$-function of the property (D)subshift $\widetilde{X} \subset
\widetilde{\Sigma}^{\Bbb Z}$. Let $\varphi: X \to \widetilde{X} $ be a
bipartite coding that is given by specifications $\psi : \Sigma \hookrightarrow
\Delta \widetilde{\Delta}, \widetilde
\psi: 
\widetilde \Sigma \hookrightarrow \widetilde\Delta \Delta$. Then
$\widetilde {g}$ is the image of $g$ under $\varphi$ if and only if
$\varphi(E(g)) = E(g)$, and for  $x \in E(G), x^{-} = x_{(-\infty, 0]},
\widetilde {x}^{-} = \varphi (x)_{(-\infty, 0]},\widetilde {\sigma} \in
\Gamma^{+}_{1}(\widetilde {x}^{-} )$,
$$
\widetilde{g}(\widetilde {x}^{-},\widetilde {\sigma}) 
\sum_{\{\sigma_{1} \in \Gamma^{+}_{1}(x^{-} ): \iota^{+} \psi (\sigma_{1}) = 
\iota^{-} \widetilde{\psi}(\widetilde {x}^{-}_{0})\}} g(x^{-},\sigma_{1}) =
$$
$$
\sum_{\{ \sigma_{1}\sigma_{2 } \in \Gamma^{+}_{2}(x¡{-}): \iota^{+}
\psi(\sigma_{1}) =  \iota^{-}(\widetilde {x}^{-}_{0}), \widetilde{\sigma} =
\widetilde{\psi}^{-1}(\iota^{-}(\psi(\sigma_{1})\iota^{+}(\psi(\sigma_{2}))\}} 
g(x^{-},\sigma_{1})g((x^{-},\sigma_{1}),\sigma_{2}).
$$
\endproclaim
\demo{Proof}
Rewrite (3) for the case of a bipartite coding.
\qed
\enddemo

Call a residually defined $g$-function $g$ of a property $(D)$ subshift $X
\subset \Sigma^{\Bbb Z}$, and a residually defined $g$-function $\widetilde
{\varphi}$ of a property $(D)$ subshift $\widetilde {X} \subset 
\widetilde{\Sigma}^{\Bbb Z} $ one step strong shift equivalent if there exist
a bipartite coding  $\varphi: X \to  \widetilde {X}  $ such that (4) 
holds, introducing in this way a notion of strong shift equivalence 
for residually defined $g$-functions of property $(D)$ subshifts. As
a consequence of Nasu's theorem topological conjugacy of residually
defined $g$-functions of property $(D)$ subshifts is equivalent to strong shift
equivalence. Similarly, topological conjugacy of $g$-measures of residually
defined $g$-functions is equivalent to 
strong shift equivalence of $g$-measures, appropriately defined. For the case
of Markov measures see
here (\cite {PW}, Theorem 3.3).
\proclaim{Proposition 3.6}
Let 
$X \subset \Sigma^{\Bbb Z}, \widetilde {X} \subset 
\widetilde{\Sigma}^{\Bbb Z}, $
be property (D) subshifts and let 
$\varphi: X \to  \widetilde {X}  $
be a topological  conjuacy. Let 
$\mu$
be a g-measure of the residually defined g-function $g$ of $X$. Then 
$\varphi \mu$ is a g-measure  of $\varphi(g).$
\endproclaim
\demo{Proof}
By Nasu's theorem sufficient to 
prove the statement of the proposition for bipartite codings, and 
for this case the statement follows from Lemma 3.5.
\qed
\enddemo


\proclaim{Proposition 3.7}
Existence and unique existence of a g-measure are invariants of the 
topological conjugacy of residually defined g-functions.
\endproclaim
\demo{Proof}
Apply Propositions 3.6.
\qed
\enddemo

\heading 4. A class of compact Shannon graphs
\endheading
Let
$\Sigma$
be a finite alphabet. We set for 
$\mu, \nu  \in \Cal M (\Sigma)$
and for 
$k \in \Bbb N$
$$
d_{k}(\mu, \nu) = \max_{a \in  \Sigma ^{k}}\vert \mu  (C(a)) - \nu(C(a))\vert.
$$
A Shannon sub-graph of $\Cal M (\Sigma)$ being uniquely determined by its
vertex set, we employ for vertex sets the same terminology as for the graphs
themselves.

We formulate a condition (I) for transition complete vertex sets 
$\Cal M \subset \Cal M ( \Sigma):$
$$
\Cal M = \overline{\bigcup_ {\alpha \in \Sigma } \tau (\alpha) \Cal M }.
\tag I
$$
\proclaim {Lemma 4.1}
A transition complete vertex set 
$\Cal M \subset \Cal M ( \Sigma)$
that satisfies condition
(I)
presents a subshift.
\endproclaim
\demo{Proof}
A word
$a \in \Sigma^{L}, L\in \Bbb N,$
is the label sequence of a finite path in $\Cal M$ precisely if there is a
$\mu \in \Cal M$
such that 
$\mu (C(a)) > 0.$
By conditon (I) there is then  an 
$\alpha \in \Sigma$
and a
$\nu \in \ \Cal M$
such that 
$$
d_{L+1}(\Phi(\alpha)\nu, \mu) < \frac{1}{2}  \mu((a))^{2}.
$$
Then 
$\nu(C(\alpha, a)) > 0$
and it is seen that the language of label sequences of finite paths in $\Cal
M$  is bi-extensible, what had to be proved.
\qed
\enddemo

The subshift that is presented by a transition complete vertex set 
$\Cal M  \subset\Cal M  (\Sigma)$ that satisfies condition (I) is 
denoted by $X^{(\Cal M)}.$

\proclaim {Lemma 4.2}
Let 
$\Cal M \subset \Cal M(\Sigma)$
be a transition complete vertex set that satisfies condition (I). Then
$$
\Cal M = \overline {\bigcup_{ a \in X^{(\Cal M)}_{(I, 0]}} \tau (a)\Cal M }  
, \qquad I \in \Bbb N.
$$
\endproclaim
\demo{Proof}
The proof is by an induction that starts with condition (I). For the induction
step, assume for an 
$I \in \Bbb N$
that
$$
\Cal M = \overline {\bigcup_{ a \in X^{(\Cal M)}_{(I, 0]} }\tau (a)\Cal M }, 
\tag 2
$$
and let 
$\mu \in \Cal M, k \in \Bbb N, \epsilon > 0.$
By  (2) there is then an 
$ a \in X^{(\Cal M)}_{(I, 0]} $
and a
$\nu \in \Cal M$
such that
$$
d_{k}(\Phi(a)\nu, \mu) < \epsilon, \tag 3
$$
and by condition (I) there is then a 
$\kappa \in  \Cal M $
and an 
$\alpha \in \Sigma$
such that 
$$
d_{k+1}(\Phi (\alpha) \kappa, \nu) < \kappa(C(\alpha))^{-1} \epsilon. \tag 4
$$
(3) and (4) imply that 
$$
d_{k+1}(\tau (\alpha a) \kappa, \nu) <  \epsilon.\qed
$$ 
\enddemo

 
Given a transition complete vertex set 
$\Cal M  \subset\Cal M  (\Sigma)$ that satisfies condition (I) we set 
for 
$x^{-} \in X^{(\Cal M)}  _{(- \infty, 0]}$,
$$
\Cal M (x^{-}) = \bigcap _{I \in \Bbb N}  \overline{\Phi(x^{-}_{(-I, 
0]})    \Cal M }.
$$
\proclaim {Lemma 4.3}
Let
$\Cal M \subset \Cal M (\Sigma)$
be a transition complete vertex set that satisfies condition (I). Then
$$
\Cal M = \bigcup_{x^{-} \in X^{(\Cal M)}_{( - \infty, 0 ]} }\Cal M (x^{-}).
$$
\endproclaim
\demo{Proof} Let 
$\mu \in \Cal M$.
Denote for $k\in \Bbb N$ by 
$\Cal A^{(k)}$ the set of 
$a \in X^{(\Cal M)}_{(- k, 0]}$
such that there is a
$\nu \in \tau (a) \Cal M$
such that
$$
d_{k}(\nu, \mu) < 2^{-k}.
$$
By Lemma 4.2 ,
$$
\Cal A^{(k)} \neq \emptyset,\qquad  k \in \Bbb N.
$$
Also
$$
\Cal A^{(k+l)}   \subset \Cal A^{(k)}, \qquad  k,l \in \Bbb N.
$$
Set
$$
\Cal B^{(k)} = \bigcap_{l \in \Bbb N} \Cal A^{(k+l)}_{( - k, 0 ]}  , \qquad k
\in \Bbb N.
$$
One can determine an
$
x^{-} \in X^{(\Cal M)}_{(- \infty, 0]}, \ x^{-}  = (\beta_{i}  )_{- \infty \leq
i
\leq 0},
$
by choosing 
$\beta_{0} \in \Cal B^{(0)}$
and by choosing inductively
$
\beta_{-k} \in \Sigma , k \in \Bbb N,
$
such that
$$
(\beta _{i})_{-k \leq i \leq 0} \in \Cal B^{(k)}, \qquad k \in \Bbb Z_{+}.
$$
Then
$x^{-}\in X^{[\Cal M)}_{(- \infty, 0]},$
and there are
$$
\nu^{(k)} \in \Phi((\beta_{i})_{ -k \leq i \leq 0}) \Cal M, \qquad k \in \Bbb N,
$$
such that
$$
d_{k}( \nu^{(k)}  , \mu) < 2^{-k}, \qquad k \in \Bbb N, 
$$
which proves that 
$\mu \in \Cal M(x^{-}).$
\qed
\enddemo

Given a transition complete vertex set
$\Cal M \subset \Cal M (\Sigma)$ that satisfies condition (I) set
$$
D^{-}_{n}(  \Cal M ) =\{ x^{-} \in X^{ (\Cal M )} _{(- \infty, 
0]}:  \vert \Cal M(x^{-})_{[1, n]} \vert = 1 \},
$$
and
$$
D^{-}_{\infty}(  \Cal M ) = \bigcap_{n\in \Bbb N} D_{n}( \Cal M),
$$
and for  $x^{-} \in D^{-}_{\infty}(  \Cal M )$ let $\mu(x^{-}) \in   \Cal M
(\Sigma)$  be the unique element in 
$
\Cal M (x^{-}). 
$
We say that a transition complete vertex set
$\Cal M \subset \Cal M (\Sigma)$
is residually contractive if it satisfies conditon (I) and if it satisfies the
following
conditions (II) and (III):
$$
\overline {D^{-}_{\infty}  (  \Cal M)} = X^{ (\Cal M )} _{(- \infty, 
0]}. \tag II
$$
\noindent
(III) For $x^{-} \in X^{ (\Cal M )} _{(- \infty, 
0]} $ and $\mu \in \Cal M(x^{-})$ there exist
$x^{-}(k) \in  D_{\infty}(  \Cal M )  , k \in     \Bbb N$, such that 
$$
\lim_{k \to \infty}x^{-}(k) = x^{-}, \quad \lim_{k \to \infty}\mu^{(\Cal
M)}(x^{-}(k)) = \mu.
$$
\proclaim {Proposition 4.4}
For a residually defined g-function $g$ of a property (D) subshift 
$X \subset \Sigma^{\Bbb Z}$,  $\Cal M (g)$
is residually contractive and presents $X$. One has
$$
D^{-}_{k}(g) =  D^{-}_{k}(\Cal M(g)), \qquad k \in \Bbb N,  \tag 1
$$
and
$$
\Cal M(g)(x^{-}) = \{\mu_{\infty}(g)(x^{-}): x^{-} \in D^{-}_{\infty}(g)\}. \tag
2
$$ 
\endproclaim
\demo{Proof} To see that 
$\Cal M (g)$ is transition complete, let 
$\mu \in \Cal M(g),  \alpha \in \Sigma, \nu \in \Cal  M (\Sigma)$,
be such that
$$
\nu = \tau ( \alpha) \mu,
$$
and let 
$k \in \Bbb N $. $\epsilon > 0.$
To find a
$y^{-} \in D_{\infty} $ such that
$$
d_{k} (\nu, \mu(y^{-}))  < \epsilon, \tag 1
$$
note that 
$\mu((\alpha)) > 0,$ and that there is an 
$x^{-} \in D_{\infty}  $ such that
$$
d_{k+1}(\mu, \mu(x^{-})) < \mu(C(\alpha)) \epsilon. \tag 2
$$
For 
$ y^{-}  = (x^{-} ,\alpha),$ (2) implies (1).

$\Cal M(g)$
satisfies condition (I) and one checks that 
$ \Cal M(g)$
and 
$\{\mu(x^{-}): x^{-} \in D_{\infty}(g)\}$ both
present X. (1) and (2) follow
from the continuity of the mapping
$$
 x^{-}  \to \mu _{\infty}(x^{-}) \quad (x^{-} \in D_{\infty}(g)).
$$

We prove that 
$\Cal  M(g)$
satisfies condition (III). For this let  
$x^{-}\in X_{(-\infty,0]}$, $\mu \in \Cal M (g)(x^{-}),$
and let 
$\nu^{(k)} \in \Cal M(g)$
be such that
$$
d_{k}(\mu, \tau (x^{-}_{(-k,0]})\nu^{(k)} < 2^{-k}, \qquad k \in \Bbb N.
$$
Then choose 
$y^{-}(k) \in D_{\infty}(g)$
such that 
$$
d_{2k}(\nu^{(k)}, \mu(y^{-}(k) )) < \nu^{(k)}(C(x^{-}_{(-k,0]} )) 2^{-k}  .
$$
Setting
$$
x^{-}(k)  = (y^{-}(k) , x^{-}_{(-k,0]}), \qquad k \in \Bbb N. \tag 3
$$
one has 
$x^{-}(k) \in D_{\infty}(g), k \in \Bbb N,$
and by (1) and (2)
$$
d_{k}(\mu(x^{-}(k)) ,\mu) < 2^{-k},\qquad k \in \Bbb N,
$$
and by (3)
$$
\lim_{k \to \infty}  x^{-}(k)  =  x^{-} .\qed
$$

From a residually contractive vertex set $\Cal M \subset \Cal M (\Sigma)$, one
obtaines a residually defined $g$-function $g_{\Cal M}$ of the subshift
$X^{(\Cal M)}$ with $C^{-}(g_{\Cal M}) = D
_{1}(\Cal M)$, by setting 
$$
g_{\Cal M}(x^{-}, \alpha) = \mu(x^{-})(C( \alpha)), \qquad  x^{-}\in D
_{1}(\Cal M),
\alpha   \in \Sigma.
$$
In particular, for a residually contractive vertex set $\Cal M \subset \Cal M
(\Sigma)$, the subshift $X^{(\Cal M )}$ has property $(D)$.
\proclaim{Proposition 4.5}
(a) For a residually definded $g$-function $g$ of a property $(D)$ subshift $X$ 

$$
g_{\Cal M(g)} = g.  
$$
(b) For a residually contractive vertex set 
$\Cal M \subset  \Cal M(\Sigma)$,
$$
\Cal M (g_{\Cal M(g)}) = \Cal M . 
$$
\endproclaim
\demo{Proof}
This is confirmed by inspection. For (b) use Lemma 4.3 in conjunction with
conditions (II) and (III).
\qed
\enddemo
\proclaim{Corollary 4.6}
A subshift has property (D) if and only if it admits a presentation by a
residually 
contractive Shannon graph.
\endproclaim
\demo{Proof}
This follows from Proposition 3.1 and Propositions 4.5.
\qed
\enddemo
Residually contractive vertex sets can be
replaced  by equivalent structures, which are measure $\lambda$-graph systems as
described in \cite {Kr2}. The measure $\lambda$-graph systems that in this way 
arise from residually contractive vertex sets can then be appropriately
referred to as residually
contractive measure $\lambda$-graph systems. The notion of strong shift
equivalence of residually defined $g$-functions will then, via a notion of
strong shift equivalence of residually contractive vertex sets, give rise
to a notion
of strong shift equivalence of residually contractive measure $\lambda$-graph
systems.

\bigskip
{\it e-mail}: 

krieger{\@}math.uni-heidelberg.de

\bye